%% file: main.tex
\documentclass[11pt]{article}
\usepackage[english]{style}
\usepackage{latexsym,amsmath,amssymb,amsfonts,diagrams,graphics}
\input{defs}

\input{title}

\begin{document}

\maketitle

\input{abstract}
\input{intro}
\input{if}
\input{iff}

\input{biblio}
\begin{center}
\begin{small}
\begin{tabular}{ll}
Dima Arinkin & Andrei C\u ald\u araru \\
Department of Mathematics& Mathematics Department\\
University of North Carolina--Chapel Hill & University of
Wisconsin--Madison\\
CB \#3250, Phillips Hall & 480 Lincoln Drive\\
Chapel Hill, NC 27599& Madison, WI 53706\\
{\em e-mail: }{\tt arinkin@math.unc.edu} &  
{\em e-mail: }{\tt andreic@math.wisc.edu}
\end{tabular}
\end{small}
\end{center}

\end{document}

%% file: defs.tex
\DeclareFontFamily{U}{rsf}{}
\DeclareFontShape{U}{rsf}{m}{n}{
  <5> <6> rsfs5 <7> <8> <9> rsfs7 <10-> rsfs10}{}
\DeclareMathAlphabet{\mathscr}{U}{rsf}{m}{n}

\DeclareMathAlphabet{\mathgth}{U}{euf}{m}{n}

\DeclareFontFamily{U}{cyr}{}
\DeclareFontShape{U}{cyr}{m}{n}{
  <5> wncyr5 <6> wncyr6 <7> wncyr7 <8> wncyr8 <9> wncyr9 <10-> wncyr10}{}
\DeclareMathAlphabet{\mathcyr}{U}{cyr}{m}{n}

\input cyracc.def

\makeatletter
\def\operator@font{\sf}
\makeatother

\newcommand{\gog}{\mathgth{g}}
\newcommand{\gh}{\mathgth{h}}
\newcommand{\gn}{\mathgth{n}}
\newcommand{\gm}{\mathgth{m}}

\newcommand{\Ug}{{\mathbf U}\gog}
\newcommand{\Sg}{{\mathbf S}\gog}
\newcommand{\Sn}{{\mathbf S}\gn}

\newcommand{\cG}{{\mathscr G}}

\newcommand{\cI}{{\mathscr I}}

\newcommand{\cL}{{\mathscr L}}

\newcommand{\cO}{{\mathscr O}}

\newcommand{\bF}{{\overline{F}}}

\newcommand{\bV}{{\overline{V}}}

\newcommand{\bN}{{\overline{N}}}

\newcommand{\bs}{{\bar{s}}}

\newcommand{\HKR}{{\mathsf{HKR}}}

\newcommand{\sHom}{\underline{\mathsf{Hom}}}

\newcommand{\D}{{\mathbf D}}

\newcommand{\chk}{{\scriptscriptstyle\vee}}

\newcommand{\At}{{\mathsf{at}}}

\newcommand{\pt}{{\mathsf{pt}}}

\newcommand{\bbk}{\mathbf{k}}
\newcommand{\bbL}{\mathbb{L}}
\newcommand{\blob}{{\scriptscriptstyle\bullet}}
\newcommand{\Kw}{K^{\scriptscriptstyle{\wedge}}_\blob}
\newcommand{\Ko}{K^{\scriptscriptstyle{\otimes}}_\blob}

\renewcommand{\S}{{\mathbb{S}}}

\newcommand{\Hom}{{\mathsf{Hom}}}

\DeclareMathOperator{\ord}{ord}

\DeclareMathOperator{\id}{id}

\DeclareMathOperator{\chr}{char}

\DeclareMathOperator{\Ext}{Ext}

\DeclareMathOperator{\codim}{codim}
\DeclareMathOperator{\sym}{Sym}

\newcommand{\Tc}{{\mathbb{T}^c}}
\newcommand{\T}{\mathbb{T}}

\newcommand{\adjoint}{\dashv}

\newcommand{\ra}{\rightarrow}

\newcommand{\lra}{\longrightarrow}
\newcommand{\contract}{\lrcorner}

\newcommand{\Z}{\mathbf{Z}}

\newcommand{\iso}{\cong}

\newcommand{\del}{\partial}

\newcommand{\smwedge}{{\scriptstyle{\wedge}}}
\newcommand{\EG}{\mathcal{EG}}
\newcommand{\bwed}{\bigwedge\nolimits}	
\newcommand\p[1]{{{\mathbb P}^{#1}}} 
\renewcommand{\phi}{\varphi}

\newarrow{Equal} =====

%% file: title.tex
\author{%
Dima Arinkin, Andrei C\u ald\u araru}

\title{When is the self-intersection of a subvariety a fibration?}

\date{}

%% file: abstract.tex
\begin{abstract}
  We provide a necessary and sufficient condition for the derived
  self-intersection of a smooth subscheme inside a smooth scheme to be
  a fibration over the subscheme.  As a consequence we deduce a
  generalized HKR isomorphism.  We also investigate the relationship
  of our result to path spaces in homotopy theory, Buchweitz-Flenner
  formality in algebraic geometry, and draw parallels with similar
  results in Lie theory and symplectic geometry.
\end{abstract}

%% file: intro.tex
\section*{Introduction}

\paragraph
Let $Y$ be a smooth scheme and let $X$ be a smooth subscheme of $Y$.
The main goal of this paper is to give a complete answer to the
following question.
\begin{quote} 
  {\em Under what circumstances is the derived self-intersection
    $X\times_Y^R X$ a fibration over $X$?}
\end{quote}

\paragraph 
Denote the closed embedding $X \hookrightarrow Y$ by $i$ and let $W$
be the derived intersection $X\times_Y^R X$.  In this context we can
think of $W$ as the ordinary scheme $X$ with structure sheaf replaced
by a {\em structure complex} $\cO_W$, which is a commutative
differential graded algebra with certain properties.  Simplifying
things slightly for ease of exposition in the introduction we can
think of $\cO_W$ as the commutative algebra object $i^*i_* \cO_X$ of
$\D(X)$.  (All our functors are implicitly derived.)  Here $\D(X)$
denotes the derived category of coherent sheaves on $X$.

\paragraph
Any object $E\in \D(X)$ concentrated in strictly positive degrees can
be regarded as a linear fibration over $X$ by considering the
dg-scheme over $X$ with structure complex $\S(E^\chk)$, the symmetric
algebra of the dual of $E$.  Ignoring the difference between dg- and
derived, the main question of this paper can be rephrased as
\begin{quote} 
  {\em Under what circumstances is the object $i^* i_*  \cO_X$ of
    $\D(X)$ of the form $\S(E^\chk)$ for some $E\in \D(X)$?}
\end{quote}

\paragraph
Our question is motivated by the case where the inclusion
$X\hookrightarrow Y$ is the diagonal embedding $\Delta:X
\hookrightarrow X \times X$.  In this case, using an analogy with
homotopy theory, the self-intersection $W$ can be understood as the
free loop space $\cL X$.  Then $\cL X$ fibers over $X$ via the map
which associates to a loop its base point.  Equivalently one has the
HKR isomorphism~\cite{Swa},~\cite{Yek}
\[ \Delta^* \cO_\Delta \iso \S(\Omega_X^1[1]) \] 
which identifies $W$ with the linear fibration $T_X[-1]$ over $X$,
where $T_X$ is the tangent bundle of $X$.

\paragraph
Our criterion for fibering $W$ over $X$ can be expressed in two
equivalent ways, one as an extension property for the normal bundle $N
= N_{X/Y}$, and another as the vanishing of a cohomological extension
class.  We shall say the closed embedding $X\hookrightarrow Y$
satisfies condition $(*)$ if it satisfies either of the following two
equivalent properties:
\begin{center}
\begin{tabular}{cl}
& The normal bundle $N$ extends to a
vector bundle $\bN$ on the\\
$(*)$\quad\quad &  first infinitesimal neighborhood $X'$ of $X$ in $Y$.\medskip \\
& The morphism $\alpha_N$ defined below vanishes.
\end{tabular}
\end{center}

\paragraph
\label{subsec:cohostar}
For any object $V\in\D(X)$ its Atiyah class is a morphism $\At_V:
V\ra V\otimes \Omega_X^1[1]$ in $\D(X)$.  The short exact sequence
\[ 0 \ra N^\chk \ra \Omega_Y^1|_X \ra \Omega_X^1 \ra 0 \] 
gives rise to an extension class $\eta:\Omega_X^1\ra N^\chk[1]$.  The
morphism $\alpha_V$ is defined as the composition
\[
\begin{diagram}[width=2em,labelstyle=\scriptstyle,alignlabels]
\alpha_V: \quad & V & \rTo^{\At_V} &  V \otimes \Omega_X[1] &
\rTo^{\id_V \otimes  \eta[1]} &  V \otimes N^\chk[2].
\end{diagram}
\]
If $V$ is a vector bundle on $X$ the class $\alpha_V$ vanishes if and
only if there exists a vector bundle $\bV$ on $X'$ whose restriction to $X$
is isomorphic to $V$, see~\cite{HuyTho,Man}.  Thus the two conditions
in $(*)$ are equivalent.  We review this equivalence
in~(\ref{pp:threestep}).  \medskip

\noindent
The main result of this paper is the following theorem.

\begin{Theorem}
  \label{thm:mainthm}
  The derived self-intersection $W = X\times_Y^R X$ fibers over $X$ if
  and only if the closed embedding $i:X\hookrightarrow Y$ satisfies
  condition $(*)$.  In this case
  \[ \cO_W = i^* i_* \cO_X \iso \S(N^\chk[1]) = \bigoplus_j \wedge^j
  N^\chk[j], \] 
  where $N^\chk$ denotes the conormal bundle of $X$ in $Y$, hence $W$
  is identified with the linear fibration $N[-1]$.
\end{Theorem}
\medskip

\noindent
As an immediate consequence of the above theorem we obtain the
following corollary.

\begin{Corollary}
\label{cor:maincor}
  Assume the closed embedding $i:X\hookrightarrow Y$ satisfies $(*)$.
  Then the hypercohomology spectral sequence
\[ {}^2 E^{pq} = H^p(X, \wedge^q N) \Rightarrow \Hom_{\D(X)}^{p+q}(i^*
i_* \cO_X, \cO_X) = \Ext^{p+q}_Y(\cO_X, \cO_X) \]
degenerates at the ${}^2E$-page, giving rise to a direct sum
decomposition
\[ \Ext^{n}_Y(\cO_X, \cO_X) \iso \bigoplus_{p+q=n} H^p(X, \wedge^q
N). \]
\end{Corollary}
\vspace*{-6mm}

\paragraph
There are many instances where condition $(*)$ is satisfied.  This is
the case, for example, when there is a finite group $G$ acting on $Y$
and $X$ is the fixed locus of the action.  Then the normal bundle
sequence splits and its class $\eta$ vanishes, so condition $(*)$ is
trivially satisfied.  In particular this happens for the case of the
diagonal embedding $\Delta:X\hookrightarrow X \times X$ and we recover
the Hochschild-Kostant-Rosenberg isomorphism.  While this already
shows that Theorem~\ref{thm:mainthm} can be regarded as a
generalization of the HKR isomorphism, the connection between these
two results is deeper, in a sense that we make precise now.

\paragraph
By analogy with homotopy theory define spaces $\Pi_k(Y,X)$ along
with maps $X \ra \Pi_k(Y,X)$ by setting $\Pi_0(Y, X) = Y$ and
recursively for $k\geq 0$ by
\[ \Pi_{k+1}(Y, X) = X \times_{\Pi_k(Y,X)}^R X. \] 
The map $X \ra \Pi_k(Y,X)$ is the original inclusion for $k=0$, and
the diagonal for $k\geq 1$.  (We assume that we work in a context
where derived fiber products exist, be it simplicial schemes,
dg-schemes of Ciocan-Fontanine and Kapranov~\cite{CioKap}, or derived
schemes of Lurie~\cite{Lur} and To\"en--Vezzosi~\cite{ToeVez}.)

If instead of algebraic varieties $X$ and $Y$ were topological spaces
and the derived fiber product was understood as the homotopy fiber
product then the spaces $\Pi_k(Y,X)$ would be obtained through the
familiar path construction and would be identified with
\[ \Pi_k(Y, X) = \{ f: D^k \ra Y~|~f(\del D^k) \subset X \}, \]  
where $D^k$ is the $k$-ball.  These are the spaces that appear in the
definition of relative homotopy groups.

\paragraph
In the context of algebraic geometry Buchweitz and
Flenner~\cite{BucFle} proved the strongest possible version of the HKR
theorem.  If $X \ra Y$ is any morphism of schemes and if $\Delta:X
\ra X\times_Y^R X$ denotes the diagonal morphism then they construct
an isomorphism
\[ \Delta^* \cO_\Delta \iso \S(\bbL[1]), \] 
where $\bbL$ is the relative cotangent complex of the morphism $X\ra
Y$.  The usual HKR isomorphism for a smooth variety $X$ is recovered
by considering the structure morphism $X \ra \pt$, whereby $\bbL =
\Omega_X^1[0]$.

By contrast, if the map under consideration is a closed embedding of
smooth schemes, its cotangent complex is $\bbL = N^\chk[1]$.
Theorem~\ref{thm:mainthm} can be rephrased as saying that 
$i^* i_* \cO_X \iso \S(\bbL)$ if and only if $(*)$ holds.

\paragraph
The two isomorphisms of structure complexes described above can be
regarded as isomorphisms of spaces
\[ 
\begin{array}{lcll}
\Pi_1(Y,X) & \iso & \bbL^\chk & \mbox{ if condition $(*)$ holds;} \\
\Pi_2(Y,X) & \iso & \bbL^\chk[-1] & \mbox{ with no assumptions.}
\end{array}
\]
The fact that in one case we get $\bbL^\chk$ and its shift in the
other is no surprise.  In the world of algebraic topology there is a
homotopy equivalence
\[ \Omega\Pi_{i-1}(Y,X) = \Pi_i(Y,X) \] 
which is the natural analogue of the fact that $\bbL^\chk$ and
$\bbL^\chk[-1]$ are shifts of one another.  What is much more
surprising is the fact that while we need condition $(*)$ for the
space $\Pi_1(Y,X)$ to fiber over $X$, there is no condition needed for
$\Pi_2(Y,X)$.

\paragraph 
It is interesting to relate Theorem~\ref{thm:mainthm} to other similar
results.  It was noted by M.\ Musta\c t\u a that condition $(*)$ is
very similar to the condition that appears in Deligne-Illusie's
theorem~\cite{DelIll} on the degeneration of the Hodge-to-de Rham
spectral sequence.  In both cases the condition needed is the
existence of a lifting to the first infinitesimal neighborhood.  It
would be very interesting to explore this connection further.

Secondly, in symplectic geometry it was noted that the spectral
sequence that relates the singular and Floer cohomologies of a
Lagrangian submanifold which is the fixed locus of an anti-symplectic
involution is expected to degenerate.  This appears to be mirror to
the fact that condition $(*)$ holds when $X$ is the fixed locus
of an involution on $Y$ and then the spectral sequence of
Corollary~\ref{cor:maincor} degenerates.

\paragraph
Finally there is a strong analogy between the world of varieties and
that of Lie algebras, initially discovered by Kontsevich and
Kapranov.  In this theory we associate to the variety $X$ the Lie
algebra object $T_X[-1]$ of $\D(X)$, where $T_X$ is the tangent bundle
of $X$.  The analogue of Theorem~\ref{thm:mainthm} is the following
result, which will appear in~\cite{CalCalTu}.

\paragraph
{\bf Theorem (Calaque-C\u ald\u araru-Tu).\ }{\it
Let $\gh \hookrightarrow \gog$ be an inclusion of Lie algebras and
let $\gn$ denote the quotient $\gog/\gh$ as an $\gh$-module.  Then
there is a Poincar\'e-Birkhoff-Witt-type isomorphism of Lie modules
\[ \Ug/\gh\Ug \iso \Sg/\gh\Sg \iso \Sn \]
if and only if a certain condition $(**)$ is satisfied.  This
condition is the exact analogue of condition $(*)$, expressed either
as an extension property or as the vanishing of a natural cohomology class.}

\paragraph
The body of the paper consists of two parts.  In the first part we
prove the ``if'' part of Theorem~\ref{thm:mainthm} and we discuss
briefly the dependence of the resulting isomorphism on the choice of
lifting $\bN$ of the normal bundle. The second part is devoted to
proving the reverse implication and to giving an explicit example of a
closed embedding $X\hookrightarrow Y$ where condition $(*)$ is not
fulfilled.  We also discuss the relationship of our class $\alpha_N$
with a natural $L_\infty$-coalgebra structure arising on the cotangent
complex. 

\paragraph
Throughout the paper we work in slightly greater generality than that
stated in this introduction.  The subscheme $X$ is only assumed to be
a local complete intersection, and is not required to be smooth.  All
schemes are assumed to be over a field $\bbk$ of characteristic zero
or greater than $\codim X$.

\paragraph {\bf Acknowledgments.} This paper grew out of conversations
with Tony Pantev, Lev Rozansky, and Craig Westerland.  The first
author is a Sloan Research Fellow and he is grateful to the Alfred
P.~Sloan Foundation for the support.  The second author is supported
by the National Science Foundation under Grant No.\ DMS-0901224.

%% file: if.tex
\section{The first implication}

In this section we prove the ``if'' part of
Theorem~\ref{thm:mainthm}.  We also discuss how the resulting
isomorphism depends on the choice of lifting $\bN$ of $N$.

\paragraph
Let $i:X\hookrightarrow Y$ be an lci closed embedding over a field
$\bbk$ which is assumed to be of characteristic zero or greater than
the codimension of $i$.  We are interested in understanding the object
$i^* i_* \cO_X$ of $\D(X)$.  

The first observation is that $i^*i_* \cO_X$ is a commutative algebra
object in $\D(X)$.  Indeed the map from the derived tensor product to
the underived one gives a morphism in $\D(Y)$
\[ i_* \cO_X \otimes i_* \cO_X \ra i_* \cO_X \] 
which, together with the canonical map $\cO_Y \ra i_*\cO_X$, makes $i_*
\cO_X$ into a commutative algebra object of $\D(Y)$.  (More
abstractly one can argue that $i_* \cO_X = i_* i^* \cO_Y$ is an
algebra object by noting that $i_* i^*$ is a monad on $\D(Y)$,
see~\cite[8.6.1]{Wei}, and using the projection formula.)  Since $i^*$ is
monoidal, the object $i^*i_* \cO_X$ of $\D(X)$ inherits a commutative
algebra structure from the same structure on $i_* \cO_X$.

\paragraph
We pause for a moment to discuss the distinction between working in
the derived category or at dg-level.  While we phrase all our results
below at the derived level, the actual multiplication maps in all our
algebra objects arise as explicit associative maps on complexes.  Thus
our algebra objects can be regarded as dg-algebras rather than
algebras in the derived category.  We stick to using the derived
category notation for ease of exposition.

\paragraph
\label{subsec:algstr}
The lci assumption on $i$ implies that the conormal bundle
$N_{X/Y}^\chk$ of $X$ in $Y$ is a vector bundle on $X$.  Throughout
this section we shall denote this vector bundle by $E$ for
simplicity. 

A local computation analogous to the one in~\cite[Appendix
A]{CalKatSha} shows that the cohomology sheaves of $i^* i_* \cO_X$
are the same as those of the symmetric algebra $\S(E[1])$,
\[ \S(E[1]) = \bigoplus_k \wedge^k E[k]. \]

\begin{Theorem}
\label{thm:ieqj}
Let $X'$ denote the first infinitesimal neighborhood of $X$ in
$Y$.  Assume that there exists a vector bundle $F$ on $X'$ whose
restriction to $X$ is isomorphic to the conormal bundle $E$ of $X$ in
$Y$.  Then there exists an isomorphism $I$ of commutative algebra
objects in $\D(X)$
\[ I:i^* i_* \cO_X \iso \S(E[1]). \] 
The isomorphism $I$ may depend on the choice of the bundle $F$.
\end{Theorem}

\paragraph 
The proof of Theorem~\ref{thm:ieqj} is divided into two parts.  In the
first part, using the existence of the vector bundle $F$ we
construct a global morphism of algebra objects
\[ I:i^*i_* \cO_X \ra \S(E[1]). \] 
Determining whether a morphism in $\D(X)$ is an isomorphism is a local
question (Lemma~\ref{lem:chkiso} below).  In the second part of the
proof we describe the restriction of the morphism $I$ to a
sufficiently small open set, and we argue that this restriction is an
isomorphism.

\paragraph
The crux of the proof lies in a computation which is familiar in the
category of vector spaces, but which makes sense equally well in
$\D(X)$.  For the sake of clarity of exposition we shall present it in
the category of vector spaces and we leave it to the reader to fill
in the details for $\D(X)$.

\paragraph
\label{subsec:calcvect}
Let $V$ be a finite dimensional vector space.  The graded vector space
$\oplus_{k \geq 0} V^{\otimes k}$ can be endowed with two distinct
algebra structures, the free algebra $\T(V)$ and the shuffle
product structure on the free coalgebra $\Tc(V)$.  The product in
$\T(V)$ is the usual, non-commutative product,
\[ (v_1 |\cdots | v_p) \cdot (v_{p+1} | \cdots |
v_{p+q}) = v_1 | \cdots | v_{p+q}, \] 
while multiplication in $\Tc(V)$ is commutative, given by the shuffle
product~\cite[6.5.11]{Wei}
\[ (v_1 | \cdots | v_p) \cdot (v_{p+1} | \cdots |
v_{p+q}) = \sum_{\sigma\mathrm{\ is\  a\ }p-q-\mathrm{shuffle}} v_{\sigma_1}
| \cdots | v_{\sigma_{p+q}}. \] 
(We have used a vertical bar to denote tensor products between
elements in $V$.)

The symmetric algebra $\S(V)$ is naturally a subalgebra of $\Tc(V)$ and
a quotient algebra of $\T(V)$.  The inclusion and quotient maps are
given by
\begin{align*}
\S(V) \ra \Tc(V) & \quad v_1 v_2 \cdots v_k \mapsto \sum_{\sigma \in
  \Sigma_k} v_{\sigma_1} | v_{\sigma_2} | \cdots |
v_{\sigma_k},
\intertext{and}
\T(V) \ra \S(V) & \quad v_1 | v_2 | \cdots | v_k \mapsto
v_1v_2\cdots v_k.
\end{align*}

\paragraph
Define the exponential map to be the isomorphism $\exp: \Tc(V)
\xrightarrow{\sim} \T(V)$ which multiplies by $1/k!$ on $V^{\otimes
  k}$.  (It is here that we need to assume that $\chr(\bbk) = 0$.)
The heart of the proof of Theorem~\ref{thm:ieqj} is the observation,
which can be confirmed by a quick calculation, that the composition
\[ \Tc(V) \xrightarrow{\exp} \T(V) \ra \S(V) \]
is a commutative ring morphism splitting the natural map $\S(V) \ra
\Tc(V)$, i.e., the composition
\[ \S(V) \ra \Tc(V) \xrightarrow{\exp} \T(V) \ra \S(V) \]
is the identity.

\paragraph
We now return to the problem of understanding $i^* i_* \cO_X$.  Begin
by fixing notation.  Let $\cI$ denote the ideal sheaf of $X$ in $Y$,
and let $X'$ be the first infinitesimal neighborhood of $X$ in
$Y$, i.e., the subscheme of $Y$ cut out by $\cI^2$.  The map $i$ then
factors into the closed embeddings
\[ X \xrightarrow{f} X' \xrightarrow{g} Y. \] 
The coherent sheaf $\cI/\cI^2$, regarded as a sheaf on $X$, is the
conormal bundle $E:=N_{X/Y}^\chk$.  Fix a vector bundle $F$ on
$X'$ together with an isomorphism $F|_X \iso E$.  The exterior
powers of $E$ will be denoted by $E^k := \wedge^k E$, and similarly
for $F$.  
\medskip

\noindent
The following proposition shows that the free coalgebra arises
naturally in the context of studying $i^* i_* \cO_X$.

\begin{Proposition}
\label{prop:iprime}
A choice of bundle $F$ and isomorphism $F|_X \iso E$ determines an
isomorphism of commutative algebras
\[ f^* f_* \cO_X \xrightarrow{\sim} \Tc(E[1]), \] 
where $\Tc(E[1])$ denotes the free coalgebra on $E[1]$, i.e., the object
\[ \Tc(E[1]) = \bigoplus_{k\geq 0} E^{\otimes k}[k] \]
of $\D(X)$, endowed with the shuffle product.
\end{Proposition}

\begin{Proof}
In order to compute $f^* f_* \cO_X$ we need to construct an explicit
resolution $\Ko$ of $f_*\cO_X$ on $X'$.  Regarding the short
exact sequence of sheaves on $Y$
\[ 0 \ra \cI/\cI^2 \ra \cO_Y/\cI^2 \ra \cO_Y/\cI \ra 0 \]
as a sequence of sheaves on $X'$ we get
\[ 0 \ra f_* E \ra \cO_{X'} \ra f_*\cO_X \ra 0. \]
Consider the composite morphism of sheaves on $X'$
\[ F \ra f_* (F|_X) \iso f_* E \ra \cO_{X'},  \]
whose cokernel is $f_*\cO_X$.  Its dual is a section $s$ of
$F^\chk$ which vanishes precisely along $X$, and the morphism above
can be regarded as the operation $\contract s$ of contracting $s$.

The complex $\Ko$ of sheaves on $X'$ can be written as
follows
\[ \Ko:\quad \cdots \xrightarrow{\contract s} F^{\otimes k}
\xrightarrow{\contract s} F^{\otimes (k-1)} \xrightarrow{\contract
  s} \cdots \xrightarrow{\contract s} \cO_{X'} \ra 0, \]
where the maps $\contract s$ are given by
\[ \contract s: F^{\otimes k} \ra F^{\otimes (k-1)} \qquad
(\contract s)(f_1 | \cdots | f_k) = s(f_1)
\cdot f_2| \cdots | f_k. \] 
Since the kernel of the original map $\contract s:F \ra
\cO_{X'}$ is $f_* E$, the kernel of the differential at the
$k$-th step in the above complex is immediately seen to be
\[ f_* E \otimes F^{\otimes (k-1)} \iso f_* E^{\otimes k}. \] 
It is a straightforward exercise to verify that $\Ko$ is exact
everywhere except at the last step, where the cokernel is
$f_*\cO_X$.  Thus $\Ko$ is a resolution of $f_*\cO_X$ on
$X'$. 

Having a resolution of $f_* \cO_X$ on $X'$ allows us to compute
$f^* f_* \cO_X$, thus yielding an isomorphism
\[ I': f^* f_* \cO_X \iso \bigoplus_k E^{\otimes k}[k] =
\Tc(E[1]), \] 
as objects of $\D(X)$.  Indeed, since $s$ vanishes along $X$, all the
differentials in the restriction of $\Ko$ to $X$ vanish.
Since $F^{\otimes k}|_X \iso E^{\otimes k}$, we get the above
isomorphism.  To finish the proof of Proposition~\ref{prop:iprime} we
only need to check that $I'$ respects the algebra operations
on the two sides.

For $k\geq 0$ consider the map of vector bundles
\[ \star_k : \bigoplus_{p+q=k} F^{\otimes p} \otimes F^{\otimes q}
\ra F^{\otimes k}\] 
given by the formula of the shuffle product defined earlier.  It is
straightforward to check that the collection of maps $\star_\blob$
is a chain map of complexes $\Ko \otimes \Ko \ra \Ko$ which represents
the natural map $f_* \cO_X \otimes f_* \cO_X \ra f_* \cO_X$ in
$\D(X)$.  Indeed, we only need to check the equality of the two
compositions in the diagram below:
\[
\begin{diagram}[height=2em,width=2em,labelstyle=\scriptstyle,alignlabels]
(x_1 | \cdots | x_p) \otimes (x_{p+1} | \cdots
| x_{p+q}) & \rMapsto^d & s(x_1) \cdot (x_2 | \cdots| x_p)
\otimes (x_{p+1} | \cdots | x_{p+q}) \\
& & \ \ +(x_1|\cdots |x_p) \otimes s(x_{p+1})\cdot
(x_{p+2} | \cdots | x_{p+q}) \\
\dMapsto^\star & & \dMapsto^\star \\
\sum_{\sigma\mathrm{\ is\  a\ }p-q-\mathrm{shuffle}} x_{\sigma_1} | \cdots
| x_{\sigma_{p+q}} & \rMapsto^d & \sum_{\sigma\mathrm{\ is\  a\
  }p-q-\mathrm{shuffle}} s(x_{\sigma_1}) \cdot x_{\sigma_2} | \cdots
| x_{\sigma_{p+q}}. 
\end{diagram}
\]
(We have written the maps without signs for simplicity, but in fact
all the morphisms have signs in them, obtained by the rule that says
that the $x_i$'s behave like odd elements in a graded vector space.)

Having represented the multiplication map $f_* \cO_X \otimes f_* \cO_X
\ra f_*\cO_X$ by the chain map $\star_\blob$, the product map of the
algebra $f^* f_* \cO_X$ is obtained by applying the functor $f^*$ to
$\star_\blob$.  It is obvious from the definition of $\star_\blob$
that the induced map under the isomorphism
\[ f^* f_* \cO_X \iso \bigoplus_k E^{\otimes k} \]
is given precisely by the shuffle product, thus $I'$ is an algebra
map.  
\qed
\end{Proof}

\paragraph
Our candidate map for the isomorphism 
\[ I: i^* i_* \cO_X \ra \S(E[1]) \]
is the composite map
\[ i^* i_* \cO_X = f^* g^* g_* f_* \cO_X \ra f^* f_* \cO_X
\xrightarrow{\sim} \Tc(E[1]) \xrightarrow{\exp} \T(E[1]) \ra
\S(E[1]), \] 
where the first morphism is the counit of the adjunction $g^* \adjoint
g_*$, the middle one is multiplying by $1/k!$ on $E^{\otimes k}[k]$,
and the last one is the natural projection map.  It is obvious that
$I$ is an algebra map: the counit map is an algebra map by standard
facts about monads, and the composition $\Tc(E[1]) \xrightarrow{\exp}
\T(E[1]) \ra \S(E[1])$ is a homomorphism by a calculation entirely
analogous to the one for vector spaces in~(\ref{subsec:calcvect}).
Thus in order to complete the proof of Theorem~\ref{thm:ieqj} we only
need to argue that $I$ is an isomorphism in $\D(X)$.  This can be
checked locally, as the following lemma shows.

\begin{Lemma}
\label{lem:chkiso}
Let $X$ be a scheme, and let $f:A\ra B$ be a morphism of objects in
$\D(X)$.  Then $f$ is an isomorphism if and only if $f|_U:A|_U \ra
B|_U$ is an isomorphism in $\D(U)$ for all open sets $U$ in a covering
of $X$.
\end{Lemma}

\begin{Proof}
The map $f$ is an isomorphism if and only if the induced maps
$H^\blob(f)$ on cohomology sheaves $H^\blob(A) \ra H^\blob(B)$ are
isomorphisms.  Since checking that a map of sheaves is an isomorphism
is a local question, and the maps on cohomology commute with
restriction, the result follows.
\qed
\end{Proof}

\paragraph
In order to complete the proof of Theorem~\ref{thm:ieqj} we shall
compare two resolutions of $\cO_X$, one on $Y$ and one on $X'$.
If $X$ were cut globally out of $Y$ by a section $\bs$ of a vector
bundle $\bF$, then the Koszul resolution of the section $\bs$ would be
a resolution of $\cO_X$ on $Y$.  This may not be the case globally,
but the following lemma shows that such a vector bundle and section
always exist locally.

\begin{Lemma}
Let $(R, \gm)$ be a regular local ring and let $J$ be an ideal in $R$ of
height $k\geq 1$.  Let $f_1,\ldots, f_k$ be elements of $R$ whose
reduction modulo $J^2$ generate the $R/J^2$-module
$J/J^2$.  Then $f_1,\ldots, f_k$ are a regular sequence for the
$R$-module $J$.
\end{Lemma}

\begin{Proof}
  The elements $f_1\mod J^2, \ldots, f_k\mod J^2$ generate $J/J^2$,
  hence their reductions mod $\gm$ generate the $R/\gm$ vector space
  $(J/J^2)/\gm = J/\gm$.  By Nakayama's lemma, $f_1, \ldots f_k$
  generate $J$.  Since $J$ was assumed to be of height $k$, the result
  follows by Krull's Hauptidealsatz.  \qed
\end{Proof}

\paragraph
We want to show that the restrictions of the map $I$ to small enough
open sets in $X$ are quasi-isomorphisms.  The previous lemma shows
that by replacing $Y$ with a sufficiently small neighborhood around
any of its points we can assume that there exists a vector bundle
$\bF$ on $Y$ and a section $\bs$ of $\bF$ such that the restriction of
$\bF$ to $X'$ is isomorphic to $F$, and the section $\bs$ maps to $s$
under this isomorphism.  Consider the Koszul resolution
\[ \Kw:\quad 0 \ra \bF^k \xrightarrow{\contract \bs} \bF^{k-1}
\xrightarrow{\contract \bs} \cdots \xrightarrow{\contract \bs} \bF
\xrightarrow{\contract \bs} \cO_Y \ra 0, \] 
where we have denoted by $\bF^k$ the vector bundle $\wedge^k \bF$ on
$Y$.  The differentials in $\Kw$ can be written down explicitly,
in a similar way we did for $\Ko$:
\[ \contract \bs:\wedge^k \bF \ra \wedge^{k-1} \bF \qquad (\contract
\bs)(f_1\smwedge \cdots \smwedge f_k) = \sum_{i=1}^k (-1)^i \bs(v_i) v_1
\smwedge \cdots \smwedge \hat{v}_i \smwedge \cdots \smwedge v_k. 
\]

Restricting $\Kw$ to $X'$ yields a non-exact complex, 
representing $g^* g_* f_* \cO_X$.  The adjunction $g^* \adjoint
g_*$ yields the counit map
\[ g^* g_* f_* \cO_X \ra f_* \cO_X. \] 
Since $\Ko$ is a resolution of $f_* \cO_X$ on $X'$ it is
reasonable to search for a map of complexes
\[ \Kw|_{X'}  \ra \Ko. \]

A straightforward computation shows that the following
diagram of bundles on $X'$ commutes for any $k\geq 0$
\[ 
\begin{diagram}[height=2em,width=2em,labelstyle=\scriptstyle]
F^k & \rTo^{\contract \bs} & F^{k-1} \\
\dTo_{\epsilon} & & \dTo_{\epsilon} \\
F^{\otimes k} & \rTo^{\contract s} & F^{\otimes (k-1)}, 
\end{diagram}
\]
where the top map is the restriction of the map $\contract \bs$ to
$X'$, and the vertical maps are the symmetrizations maps defined
below:
\[ \epsilon(v_1\smwedge  v_2 \smwedge \cdots\smwedge v_k) = \sum_{\sigma\in \Sigma_k}
(-1)^{\epsilon(\sigma)} v_{\sigma_1} | v_{\sigma_2} | \cdots | v_{\sigma_k}.
\]
(This map is the analogue of the map $\S(V) \ra \Tc(V)$ discussed
earlier, with signs included to take into account the fact that
$F[1]$ in the symmetric monoidal category $\D(X)$ behaves like an odd
vector space in the category of $\Z_2$-graded vector spaces.)  Thus we
get a map of complexes $\epsilon_\blob$ from the restriction of the
Koszul resolution to $X'$ to the second resolution, and it is
immediate to check that it represents the counit map above.

Restricting further the map $\epsilon_\blob:\Kw \ra \Ko$
to $X$ by pulling-back via $f$, all the differentials in both
complexes vanish.  Therefore we conclude that the natural map
\[ i^* i_* \cO_X \ra f^* f_* \cO_X \] 
arising by contracting $g^* g_*$ is given by the
map of complexes $\epsilon : \S(E[1]) \ra \Tc(E[1])$.  Thus composing
the map $I:i^* i_* \cO_X \ra \S(E[1])$ with the (local) isomorphism
$\S(E[1]) \xrightarrow{\sim} i^* i_* \cO_X$ obtained from the restriction
of the Koszul resolution gives the composition
\[ \S(E[1]) \ra \Tc(E[1]) \xrightarrow{\exp} \T(E[1]) \ra \S(E[1]) \]
which we know is an isomorphism.  Therefore the maps induced by $I$
are locally isomorphisms on cohomology sheaves, and thus $I$ is globally a
quasi-isomorphism.  This concludes the proof of
Theorem~\ref{thm:ieqj}.  
\qed

\paragraph
There are two clear situations in which we can guarantee the existence
of a vector bundle $F$ on the first infinitesimal neighborhood $X'$
extending the conormal bundle $E=N_{X/Y}^\chk$.  One is when $X$ is
a global complete intersection in $Y$, that is there exists a vector
bundle $V$ on $Y$ of rank $\codim(X/Y)$ and a section of $V$ whose
vanishing locus is $X$ scheme-theoretically.  In this case $V^\chk|_X
= N^\chk$, and hence $F = V^\chk|_{X'}$ extends $E$.

\paragraph
A more interesting situation is when the sequence 
\[ 0 \ra N^\chk \ra \Omega_Y^1|_X \ra \Omega_X \ra 0 \] 
is split.  Then the class $\eta$ which appears in the definition of
the obstruction $\alpha_V$ in~(\ref{subsec:cohostar}) vanishes, and
hence any object of $\D(X)$ extends to $X'$.  This can be understood
geometrically by noting that if $\eta$ vanishes then there is a map
$\pi:X' \ra X$ splitting the natural inclusion $X\hookrightarrow X'$,
see~\cite[20.5.12 (iv)]{EGA4.1}.  Then any $V\in \D(X)$ extends to the
object $\pi^*V \in \D(X')$ whose restriction to $X$ is $V$.

\paragraph
The above sequence is split in two situations of interest.  One is
when $i:X\hookrightarrow Y$ is split by a morphism $\pi:Y\ra X$.  Such
is the case for the diagonal embedding $X\ra X\times X$ which is split
by either projection on the factors, or for the zero section of a
bundle map $E \ra X$.  Another situation where the tangent sequence
splits is when there is a finite group $G$ of order not divisible by
$\chr \bbk$ which acts on $Y$, and $X$ is the fixed locus $Y^G$.  Then
the dual map
\[ T_X \ra T_Y|_X \]
is split by the map that takes a tangent vector $v$ to $Y$ at a point
of $X$ to 
\[ \frac{1}{\ord(G)} \sum_{g\in G} g\cdot v. \]
Since the diagonal is the fixed locus of the involution which
interchanges the two factors, this gives rise to another HKR-type
isomorphism through the corresponding splitting of the tangent sequence.

\paragraph
An interesting question that arises in the study of the diagonal
embedding is whether the two HKR isomorphisms 
\[ \HKR_1, \HKR_2: \Delta^*\Delta_* \cO_X \iso \S(\Omega_X^1[1]) \]
induced by the splittings $\pi_1$, $\pi_2$ of $\Delta$ are different.
While we are unable to answer this question at the moment, we give
strong evidence that these splittings induce the {\em same} HKR
morphism.

\paragraph
By adjunction the two isomorphisms $\HKR_1$, $\HKR_2$ can be regarded
as morphisms
\[ \cO_\Delta \ra \Delta_* \S(\Omega_X^1[1]), \] 
which in turn can be regarded as natural transformations $f$ and $g$
between the identity functor of $\D(X)$ and the functor $-\, \otimes
\S(\Omega_X^1[1])$.  Instead of arguing that the maps on the level of
kernels agree, we'll prove the weaker statement that the induced natural
transformations are the same.

\paragraph
The natural transformations $f$ and $g$ give, for $E\in \D(X)$, maps
\[  E \ra E \otimes \S(\Omega_X^1[1]) = \bigoplus_i E \otimes
\Omega^i_X[i], \] 
which, for simplicity, we denote by $f$ and $g$ as well.  

The components of these maps can easily be understood (see
also~\cite{CalMuk2} for more details).  Explicitly, if $f$ is obtained
using the second projection to split $\Delta$, then the component
$f_i: E \ra E \otimes \Omega^i_X[i]$ of $f$ is
\[ f_i = \epsilon \circ (\id_{\Omega^{\otimes(i-1)}_X[i-1]}
\,\otimes\, \At_E)\circ (\id_{\Omega^{\otimes(i-2)}_X[i-2]} \,\otimes\, \At_E)
\circ \cdots \circ \At_E, \] 
where $\epsilon$ is the antisymmetrization map $\Omega^{\otimes i}_X
\ra \Omega^i_X$.
Similarly, if $g$ corresponds to the first projection, then $g_i$ is
given by 
\[ g_i = \epsilon \circ \At_{E \otimes \Omega_X^{\otimes (i-1)}[i-1]} \circ
\At_{E \otimes \Omega_X^{\otimes (i-2)}[i-2]} \circ \cdots \circ
\At_E. \]
The fact that these two maps are the same follows easily from the fact
that the Atiyah class gives rise to a Lie coalgebra structure on
$\Omega_X^1[1]$, which makes $E$ into a Lie comodule.  

\paragraph
It is obvious that $f_i=g_i$ for $i=0,1$.  We exemplify the
calculation that $f_2 = g_2$, and leave the details that $f_i=g_i$ for
$i\geq 3$ to the reader.  Denoting by $\gog = T_X[-1]$, the dual of
$\Omega_X^1[1]$, and using $E^\chk$ instead of $E$, the duals of the
maps $f_2$ and $g_2$ are the maps
\begin{align*}
\gog \otimes \gog \otimes E & \quad \lra \quad \gog \otimes E \quad \lra
\quad E
\intertext{given by, respectively,}
x\otimes y\otimes e & \mapsto x\otimes  (y\cdot e)  \ \mapsto x\  \cdot (y
\cdot e)
\intertext{and}
x \otimes y \otimes e & \mapsto  x\cdot (y \otimes e) 
 = [x,y] \otimes e + y\otimes (x\cdot e) \\
& \quad\quad\quad\quad\quad\quad\mapsto [x,y] \cdot e  +
y\cdot (x \cdot e).
\end{align*}
(We have written the above maps in component notation for clarity, but
these maps make sense in an arbitrary symmetric monoidal category,
like $\D(X)$.)  Note the the equality 
\[ x\cdot (y \otimes e) = [x,y] \otimes e + y\otimes (x\cdot e) \] 
is nothing but the way $\gog$ acts on the tensor product of
representations $\gog\otimes E$.  In other words the above map
$\gog\otimes \gog \otimes E \ra \gog\otimes E$ is the dual of
$\At_{\Omega^1_X[1] \otimes E}$.

The equality $f_2= g_2$ follows now from the fact that $E$ is a
representation of $\gog$ in $\D(X)$.   See~\cite{RobWil} for more details.

%% file: iff.tex
\section{The reverse implication}

\newcommand{\LL}{\mathbb L}
\newcommand\W{M^\bullet}				
\newcommand\WX{M^\bullet_X}				
\newcommand\WV{M^\bullet_{(V)}}
\newcommand\RV{C^\bullet} 
\newcommand\RVX{C^\bullet_X} 

In this section we argue that if $i^* i_* \cO_X$ is formal (isomorphic
to $\S(N^*[1])$), then the obstruction class $\alpha_N$
from~(\ref{subsec:cohostar}) vanishes.  We then give an explicit
example of a closed embedding of smooth varieties for which $\alpha_N$
does not vanish.  In a final part of this section we discuss the
relationship of the class $\alpha_N$ to an $L_\infty$-coalgebra
structure on the relative cotangent complex.

\paragraph 
\label{subsec:pushpulltriangle}
As before let $i:X\hookrightarrow Y$ be a closed lci embedding with
conormal bundle $E = N^\chk_{X/Y}$ and let $V$ be a vector bundle on
$X$.  Consider the truncation $\tau^{\ge -1}(i^*i_* V)$.  It has only
two non-trivial cohomology sheaves, $H^0$ and $H^{-1}$, which are
naturally isomorphic to $V$ and $V\otimes E$, respectively.  Therefore
$\tau^{\ge -1}(i^*i_* V)$ fits into a triangle
\[ V\otimes E[1]\to\tau^{\ge-1}i^*i_*V\to
V\to V\otimes E[2]. \]

\begin{Definition} 
\label{df:alpha} 
The rightmost map $V\to V\otimes E[2]$ in the above triangle will be
called the \emph{HKR class} of $V$, denoted by
$\alpha_V\in\Ext^2(V,V\otimes E)$.
\end{Definition}

\paragraph
{\em A priori} the notation in the above definition seems ambiguous,
as we have already defined $\alpha_V$ using a different formula
in~(\ref{subsec:cohostar}).  However, the bulk of this section is
devoted to proving the following two facts: 
\begin{itemize} 
\item[a)] the HKR class $\alpha_V$ equals $\At_V \circ \eta$, the
  previously defined $\alpha_V$;
\item[b)] $\alpha_V$ is the obstruction to lifting $V$
to $X'$, the first infinitesimal neighborhood of $X$ in $Y$.
\end{itemize}

\paragraph
It is easy to see that the HKR class depends only on the embedding
$f:X\to X'$ of $X$ into its first infinitesimal neighborhood. Indeed,
the truncation of the natural map $i^*i_* V\to f^*f_*V$ is an
isomorphism
\[\tau^{\ge-1}i^*i_*V\stackrel{\sim}{\lra}\tau^{\ge-1}f^*f_*V. \]
The claim that the HKR class $\alpha_V$ is the obstruction to
extending $V$ to a vector bundle $\bV$ on $X'$ is naturally
formulated in the language of gerbes: locally on $X$ an extension of
$V$ always exists, and the local extensions form a gerbe over $X$ for the
sheaf $\sHom(V,V\otimes E)$.  The claim is that $\alpha_V$ is the class of
this gerbe.

\paragraph
Let us be more precise about these definitions. For an open subset
$U\subset X$ let $U'\subset X'$ be its first infinitesimal
neighborhood. Denote by $\EG(U)$ the category of extensions of $V|_U$
to a vector bundle on $U'$.  Thus objects of $\EG(U)$ are vector
bundles $\bV_U$ on $U'$ equipped with isomorphisms $(\bV_U)|_U\iso
V|_U$.  (To simplify notation we will usually omit this isomorphism.)
Clearly $\EG(U)$ is a groupoid. As $U$ varies the groupoids $\EG(U)$
form a sheaf of groupoids $\EG$.

If $U$ is small enough an extension $\bV_U$ exists; moreover, any two
extensions are locally isomorphic.  Finally the automorphism group of
an extension $\bV_U$ equals $\Hom(V|_U,(V\otimes E)|_U)$. In
other words $\EG$ is a gerbe over $\sHom(V,V\otimes E)$.

\begin{Proposition} 
\label{pp:directimage} 
The class
\[[\EG]\in H^2(X,\sHom(V,V\otimes E))=\Ext^2(V,V\otimes E)\] 
equals $\alpha_V$ as defined in~(\ref{df:alpha}).
\end{Proposition}

\begin{Proof} 
The statement is probably well known; since we were unable to locate a
reference, we include a detailed proof.

Let $C^\bullet =C^{-1}\stackrel{d}{\lra} C^0$ be a two-term complex of
$\cO_X$-modules with cohomology sheaves $H^0$, $H^{-1}$, such that
$H^0$ is locally free.  Then $C^\bullet$ gives rise to a gerbe $\cG$
over $\sHom(H^0, H^{-1})$ as follows.  For every open set $U$ there
is an associated groupoid $\cG(U)$.  Its objects are elements $\phi$
of $\Hom(H^0|_U,C^0|_U)$ that lift the identity automorphism of
$H^0|_U$; equivalently $\phi$ is a section of $C^0\to H^0$ over $U$.
Morphisms between between $\phi,\phi'\in \Hom(H^0|_U,C^0|_U)$ are
elements $\psi\in\Hom(H^0|_U,C^{-1}|_U)$ such that
$d\psi=\phi'-\phi$. As $U$ varies the groupoids $\cG(U)$ naturally
form a gerbe over $\sHom(H^0,H^{-1})$ whose class
\[ [\cG] \in \Hom(H^0, H^{-1}[2]) \]
is precisely the rightmost map in the triangle
\[ H^{-1}[1] \ra C^\bullet \ra H^0 \ra H^{-1}[2]. \] 
Essentially $\cG$ is the gerbe of splittings of the complex
$C^\bullet$.

Now we will pick a specific two-term complex $C^\bullet$ which will represent
$\tau^{\geq -1}(i^*i_* V)$.  Choose a truncated resolution of $i_* V$
of the form
\[0\to F^{-1}\to F^0\to i_*V\to 0,\] 
where $F^i$ are $\cO_Y$-modules (not necessarily quasi-coherent) and
$F^0$ is flat. Then the restriction
\[F^\bullet|_X=(F^{-1}|_X\to F^0|_X)\] 
represents the object $\tau^{\ge -1}(i^*i_*V)$.  Note that $F^{-1}$ is
not flat, and $F^{-1}|_X$ refers to the naive (non-derived)
restriction. Let $\cG$ be the gerbe of splittings of the complex
$F^\bullet|_X$.

It remains to construct a morphism $W:\cG\to\EG$ of $\sHom(V,V\otimes
E)$-gerbes.  To simplify the notation, we will only construct it for
global sections.  On objects, a morphism $\phi:V\to F^0$ defines an
extension $W_\phi\in\EG(X)$ by
\[W_\phi=\{s\in F^0:s|_X\in\phi(V)\}/\cI_XdF^{-1},\] 
where $\cI_X\subset \cO_Y$ is the ideal sheaf of $X$.  Given two
morphisms $\phi,\phi':V\to F^0$ and $\psi:V\to F^{-1}$ such that
$d\psi=\phi'-\phi$, the corresponding map $W_\phi\to W_{\phi'}$ sends
$s\in W_\phi$ to $s+d\psi([s])$ where $[s]$ is the image of $s\in
F^0$ in $V$. It is easy to see that we obtain a morphism of gerbes in
this way.  Thus $\alpha_V = [\cG] = [\EG]$.
\qed
\end{Proof}

\paragraph
We will also be interested in the following modification of the above
extension problem.  Let $\tilde V$ be a vector bundle (or a
quasi-coherent sheaf) on $X$, equipped with a morphism $m:V\otimes
E\to\tilde V$.  We can then consider the problem of constructing exact
sequences
\[ 0\to i_*\tilde V\to \bV\to i_*V\to 0 \]
on $Y$ such that the sheaf of ideals $\cI_X\subset \cO_Y$ acts on $\bV$ as
\[\cI_X\otimes \bV\to \cI_X\otimes i_*V=i_*(V\otimes E)\stackrel{i_*m}{\lra}
i_*\tilde V\to \bV.\] 
It is easy to see that in this case the obstruction is a class
\[\alpha_{V,m}\in \Ext^2(V,\tilde V)\] 
that is functorial in the pair $(\tilde V,m)$. Therefore
$\alpha_{V,m}=m \circ\alpha_V$.

\begin{Corollary} 
\label{co:functoriality} 
An extension 
\[ 0\to i_*\tilde V\to \bV\to i_*V\to 0 \] as above exists if and only
if $m \circ \alpha_V$ vanishes.
\end{Corollary}

\paragraph 
\label{subsec:threestep}
We would now like to relate the class $\alpha_V$ to $\eta\circ \At_V$,
which was our original definition of $\alpha_V$.  Let us assume that $X$
and $Y$ are both smooth. This assumption persists
through~\ref{subsec:endsmooth}.  It allows us to use the classical
Atiyah class instead of its derived counterpart constructed using the
cotangent complex.  See~\cite{HuyTho} for a more general statement
without this assumption.

As before $V$ is a vector bundle on $X$. Let us describe the class
$\eta \circ \At_V$ explicitly.  The Atiyah sequence of $V$ is the exact
sequence of coherent sheaves
\[ 0\to V\otimes\Omega_X^1\to J^1(V)\to V\to 0 , \] 
where $J^1(V)$ is the first jet bundle of $V$. The corresponding
extension class $\At_V\in\Ext^1(V,V\otimes\Omega_X^1)$ is called the
Atiyah class of $V$.  Composing with the class
$\eta\in\Ext^1(\Omega_X^1,E)$ of the short exact sequence
\[ 0\to E\to \Omega^1_Y|_X \to\Omega^1_X\to 0 \]
we get a morphism
\[\eta\circ\At_V \in\Hom(V, V\otimes E[2]) = \Ext^2(V,V\otimes E). \]
This composite morphism can be interpreted as the obstruction to
fitting the two exact sequences above into a diagram
\[ 
\begin{diagram}[height=2em,width=2em]
& & & & 0 & & 0 & & \\
& & & & \dTo & & \dTo & & \\
0 & \rTo & V\otimes E & \rTo & V \otimes i^* \Omega_Y^1 & \rTo &
V\otimes \Omega_X^1 & \rTo & 0 \\
& & \dEqual & & \dTo & & \dTo & & \\
0 & \rTo & V \otimes E & \rTo & J' & \rTo & J^1(V) & \rTo & 0 \\
& & & & \dTo & & \dTo & & \\
& & & & V & \rEqual & V & & \\
 & & & & \dTo & & \dTo & & \\
& & & & 0 & & 0. & & 
\end{diagram}
\]
In other words $J'$ is a vector bundle on $X$ equipped
with a filtration
\[0=J'_0\subset J'_1\subset J'_2\subset J'_3=J'\] such that
\begin{align*} J'_1&=V\otimes E&J'_2&=V\otimes
i^*\Omega^1_Y\\ J'_2/J'_1&=V\otimes\Omega_X&J'_3/J'_1&=J^1(V)\\
J'_3/J'_2&=V.
\end{align*}

\begin{Proposition} 
\label{pp:threestep} 
$\alpha_V=\eta\circ\At_V$.
\end{Proposition}

\begin{Proof} 
This is Remark~II.6.10(b) in \cite{Man} (see also Theorem~II.6.6), and
a particular case of Corollary~3.4 in~\cite{HuyTho}. Let us sketch the
argument.

Let us interpret $\alpha_V$ as the obstruction to extending $V$ to a
vector bundle $V'$ on $X'$. Such an extension $V'$ can be viewed as a
sheaf $g_*V'$ on $Y$. Considering its Atiyah extension
\[0\to g_*V'\otimes\Omega^1_Y\to J^1(g_*V')\to g^*V'\to 0 \] 
it is easy to see that the restriction $J'=i^*( J^1(g_*V'))$ fits into
the diagram of~(\ref{subsec:threestep}).  In other words, a solution to the
problem obstructed by $\alpha_V$ yields a solution to the problem
obstructed by $\eta\circ\At_V$.

Looking at local solutions, we obtain a map between the corresponding
two gerbes (the gerbes of local solutions to the two problems). Since
the gerbes are over the same sheaf $\sHom(V,V\otimes E)$,
the map is a 1-isomorphism. Therefore, the classes of the two gerbes
coincide.
\qed
\end{Proof}

\paragraph
\label{subsec:atiyahtri}
The class $\eta\circ\At_V$ can also be interpreted as
follows. Consider the composition
\[i^*\Omega_Y\otimes V\to\Omega_X\otimes V\to J^1(V),\] and let
$J^1_{X/Y}(V)$ be the length two complex
\[i^*\Omega_Y \otimes V\to J^1(V),\] 
with $J^1(V)$ located in cohomological degree zero and differential
the map above.  The cohomology sheaves of $J^1_{X/Y}(V)$ are
$V$ and $V\otimes E$ in degrees zero and minus one, respectively. We
therefore obtain a triangle
\[ V\otimes E[1]\to J^1_{X/Y}(V)\to V\to V\otimes E[2], \]
in which the rightmost map is $\eta\circ\At_V$.

Propositions~\ref{pp:directimage} and \ref{pp:threestep} can be
interpreted as saying that there is an isomorphism
$J^1_{X/Y}(V)\simeq\tau^{\ge -1}(i^*i_* V)$ that relates the triangle
above to the triangle of~(\ref{subsec:pushpulltriangle}).

\paragraph
\label{subsec:endsmooth}
\textbf{Remark. } The triangle in~(\ref{subsec:atiyahtri}) can be
viewed as a version of the Atiyah sequence if one views $X$ as a
scheme over $Y$. Here the analogue of the cotangent bundle $\omega_X^1
= \bbL_{X/\bbk}$ in the definition of the usual Atiyah class is taken
by the relative cotangent complex $E[1] = N^\chk[1] = \bbL_{X/Y}$.  We
thus see that the HKR class $\alpha_V$ is essentially the relative
Atiyah class of $V$ viewed as a sheaf on $X/Y$.

\paragraph We now want to concentrate on a special HKR class, namely
that of the conormal bundle itself.  We return to the original
assumptions that $Y$ is smooth, $i:X\hookrightarrow Y$ is a locally
complete intersection closed embedding with conormal
bundle $E$.  We are interested in the HKR class associated to $E$,
\[\alpha=\alpha_E\in H^2(X,E^{\otimes 2}\otimes E^\vee).\]

\begin{Proposition}
\label{pp:skew} 
The class $\alpha$ is skew-symmetric,
\[\alpha\in H^2(X,\wedge^2E\otimes E^\vee)\subset H^2(X,E^{\otimes
2}\otimes E^\vee).\]
\end{Proposition}
\vspace*{-6mm}

\begin{Proof}
We need to check that the image of $\alpha$ in $H^2(X,\S^2(E)\otimes
E^\vee)$ vanishes. This image equals $\sym(\alpha)$, where $\sym$ is
the symmetrization map $\sym:E^{\otimes 2}\to\S^2(E)$.  By
Corollary~\ref{co:functoriality}, we need to construct an exact
sequence
\[0\to i_*\S^2 E\to E'\to i_* E\to 0\] of sheaves on $Y$ such that
the action of $I_X$ on $E'$ is given by $\sym$. Now take
$E'=I_X/I_X^3$.
\qed
\end{Proof}

\paragraph 
We introduce the following notation. For an object $C\in\D(X)$,
consider the filtration of $C$ by objects $H^{-k}(C)[k]$. The
differential $H^{-k}(C)\to H^{-k-1}(C)[2]$ in the corresponding
spectral sequence is denoted by $\delta_k=\delta_k(C)$. Explicitly,
$\delta_k$ is the shift of the rightmost map in the exact triangle
\[H^{-k-1}(C)[k+1]\to\tau^{\le -k}\tau^{\ge -k-1} C\to H^{-k}(C)[k]\to
H^{-k-1}(C)[k+2].\] Clearly, $\delta_k(C)$ is functorial in $C$.

As an example of using this notation the HKR class of a vector bundle
$V$ on $X$ is nothing but
\[\alpha_V=\delta_0(i_*i^* V)=\delta_0(f_*f^*V).\]

\paragraph
It is obvious that if the object $C$ is formal (isomorphic to the
direct sum of its cohomology sheaves), then $\delta_k(C) = 0$ for all
$k$.  We will show that the HKR class $\alpha = \alpha_E$ of $E$ is
precisely $\delta_1(i^*i_* \cO_X)$, thus proving the second
implication of Theorem~\ref{thm:mainthm}: if $i^*i_* \cO_X$ is formal,
then 
\[ \alpha = \delta_1(i^*i_* \cO_X) = 0. \] 
Note that both $\alpha$ and $\delta_1(i^*i_*\cO_X)$ are maps from $E =
H^{-1}(i^*i_*\cO_X)$ to the shift by two of $\wedge^2E =
H^{-2}(i^*i_*\cO_X)$.

As a side note it is easy to see that the map giving the unit of the
algebra structure on $i^*i_*\cO_X$ splits the natural projection
$i^*i_* \cO_X$, so that
$i^*i_*\cO_X=\cO_X\oplus\tau^{<0}i^*i_*\cO_X$. Therefore, $\delta_0=0$
and $\delta_1$ is the first map that has a chance of not being zero.

\paragraph
As before we will understand the class $\delta_1(i^*i_* \cO_X)$ in
two steps, by first studying the behavior of this class with respect
to the embedding $f$ of $X$ into $X'$, the first infinitesimal
neighborhood of $X$ in $Y$.  Recall that $H^{-k}(f^*f_*\cO_X)=E^{\otimes k}$
and we can consider the map
\[\delta_1(f_*f^*\cO_X):E\to E^{\otimes 2}[2].\]

\begin{Lemma} 
We have 
\[ \delta_1(f_*f^*\cO_X)=\alpha; \]
here $\alpha$ is regarded as an element of $H^2(X,\wedge^2E\otimes
E^{\vee}) \subset H^2(X,E^{\otimes 2}\otimes E^{\vee})$.
\end{Lemma}

\begin{Proof} Consider the natural exact sequence
\[0\to f_*E\to\cO_{X'}\to f_*\cO_X\to 0.\] 
Applying $f^*$, we obtain an exact triangle
\[f^*f_*E\to\cO_X\to f^*f_*\cO_X\to f^*f_*E[1].\] The composition
\[\tau^{<0}f^*f_*\cO_X\to f^*f_*\cO_X \to f^*f_*E[1]\] is an isomorphism, and
therefore
\begin{equation}
\delta_1(f^*f_*\cO_X)=\delta_1(\tau^{<0}f^*f_*\cO_X)=\delta_1(f^*f_*E[1])=\delta_0(f^*f_*E)=\alpha_E.\tag*{\qed}
\end{equation}
\end{Proof}

\paragraph
Let $\epsilon:E^{\otimes 2}\to\wedge^2(E)$ be the skew-symmetrization
map. By Proposition~\ref{pp:skew}, $\alpha$ is already skew-symmetric, so that
\[\alpha=\epsilon(\alpha)\in H^2(X,\wedge^2(E)\otimes E^\vee).\]

\noindent
We are now ready to prove the main result of this section, which
finishes the proof of the second implication of Theorem~\ref{thm:mainthm}.
\medskip

\begin{Proposition} 
\label{pp:alphadelta} 
For the object $i^*i_*\cO_X\in\D(X)$ we have
\[\delta_1(i^*i_*\cO_X)=\epsilon(\alpha)=\alpha.\]
\end{Proposition}
\vspace*{-6mm}

\begin{Proof} Consider the composition
\[\sym^2(E)[2]\to E^{\otimes 2}[2]\to H^{-2}(f^*f_*\cO_X)\to \tau^{\ge
-2}(f^*f_*\cO_X),\] 
and include it into a triangle
\[\sym^2(E)[2]\to \tau^{\ge -2}(f^*f_*\cO_X)\to C\to\sym^2(E)[3].\] 
By construction, 
\[H^{-k}(C)=\begin{cases} \bwed^k E & k=0,1,2\\0
  &\text{otherwise,}\end{cases}\] 
Moreover, the map
$H^{-2}(f^*f_*\cO_X)\to H^{-2}(C)$ equals $\epsilon$, and so
\[\delta_1(C)=\epsilon(\delta_1(f^*f_*\cO_X))=\epsilon(\alpha).\]
On the other hand, the composition
\[i^* i_*\cO_X\to f^*f_*\cO_X\to \tau^{\ge -2}(f^*f_*\cO_X)\to C\]
induces an isomorphism on cohomology objects in degrees $0$, $-1$, and
$-2$. Therefore, $\delta_1(i^* i_*\cO_X)=\delta_1(C)= \alpha$. 
\qed
\end{Proof}

\paragraph
There is an alternative approach to the HKR class $\alpha$ that we
have studied above arising from a canonical Lie coalgebra structure
on the shifted cotangent complex $\LL_{X/Y}[1]$. This approach is
natural if one views $\alpha$ as a version of the Atiyah class
(Remark~\ref{subsec:endsmooth}): in the absolute case,
Kapranov~\cite{Kap} showed that the Atiyah class of a smooth manifold
$X$ gives a Lie algebra structure on the shifted tangent bundle $T_X[-1]$.

\paragraph
The coalgebra structure on the shifted cotangent complex is defined
for arbitrary morphisms of ringed spaces (or topoi). Let us sketch the
construction in this generality.

Recall the definition of the cotangent complex $\LL_{X/Y}$ for a
morphism $\phi:X \ra Y$ of ringed spaces. We follow the original
definition of Illusie~\cite{Illusie}, which also appears in the works
of Ciocan-Fontanine and Kapranov~\cite{CioKap},~\cite{Hilb}.

Consider $\cO_X$ as a sheaf of algebras over the sheaf of rings
$\phi^{-1}\cO_Y$. There exists a resolution of $\cO_X$ of the form
\[(\S(\W),d)\to\cO_X,\] 
where $\W$ is a graded flat $\phi^{-1}\cO_Y$-module concentrated in
non-positive degrees, and $d$ is a (degree one) differential on
$\S(\W)$, viewed as a graded $\phi^{-1}\cO_Y$-algebra.  Such a
resolution can be constructed by the usual iterative procedure
(`attaching cells to kill homotopy groups') as in
\cite[p.~256]{Quillen}.  This is the approach used
in~\cite[Theorem~2.6.1]{CioKap} with somewhat different
assumptions. (In~\cite{CioKap} the authors work with coherent sheaves on
quasi-projective schemes, which admit enough locally free sheaves; in
our setting there are enough flat sheaves.)  Alternatively we can
follow~\cite[I.1.5.5.6]{Illusie} and use the natural `maximal'
resolution of $\cO_X$.

Taking the tensor product we obtain a commutative dg-algebra  
\[(\S(\W)\otimes_{\phi^{-1}\cO_Y}\cO_X,d)\] 
over $\cO_X$. In the derived category of $\cO_X$-modules this algebra
represents the derived tensor product
$\cO_X\otimes_{\phi^{-1}\cO_Y}\cO_X$.

Set $\WX=\W\otimes_{\phi^{-1}\cO_Y}\cO_X$.
Let us embed $\WX$ into $\S(\W)\otimes_{\phi^{-1}\cO_Y}\cO_X$ using
the morphism
\[ \WX\to\WX\oplus\cO_X=
(\W\otimes_{\phi^{-1}\cO_Y}\cO_X)\oplus\cO_X\hookrightarrow
\S(\W)\otimes_{\phi^{-1}\cO_Y}\cO_X.\]
Here the leftmost arrow equals $(id,-\varepsilon)$, and $\varepsilon:
\WX\to\cO_X$ is the augmentation map. The embedding induces an
identification $\S(\WX)=\S(\W)\otimes_{\phi^{-1}\cO_Y}\cO_X$. It is
easy to see that the ideal
\[\S^{\ge k}(\WX)=\bigoplus_{j\ge k}\S^j(\WX)\]
is preserved by the differential $d$.  In particular $d$ induces a
differential on
\[\S^{\ge 1}(\WX)/\S^{\ge 2}(\WX)=\WX.\]
The \emph{cotangent complex} $\LL = \LL_{X/Y}$ is by definition $\WX$
equipped with this differential, viewed as an object of the derived
category of $\cO_X$-modules.

Clearly, $\WX$ depends on the choice of resolution of $\cO_X$ (unless
we use the natural resolution of~\cite{Illusie}). However, one can
check that it is unique up to quasi-isomorphism,
cf.~\cite[Corollary~II.1.2.6.3]{Illusie},~\cite[Proposition~2.7.7]{CioKap},
and~\cite[Proposition~2.1.2]{Hilb}. 

\paragraph
\label{subsec:linfty}
Recall that a square-zero differential $d$ on a free commutative
algebra $\S(V)$ is equivalent to a (possibly curved)
$L_\infty$-coalgebra structure on $V$.  (Here $V$ is a graded vector
space, or any similar object in an appropriate graded symmetric monoidal
category.)  Thus the specific complex $\WX[1]$ constructed above
inherits tautologically the structure of an $L_\infty$-coalgebra in
the dg-category of complexes of flat $\cO_X$-modules.  Its structure
maps $d_k:\WX[1]\to\wedge^k(\WX[1])[2-k]$, $k\geq 1$, are defined by
the equality 
\[d|_{\WX}=(d_1,d_2,\dots):\WX\to\S^{\ge1}(\WX)\subset\S(\WX).\] (One
can see that the curvature term $d_0:\WX \ra \cO_X$ of this coalgebra
is trivial.)

\paragraph
\textbf{Example. } 
\label{ex:bracket} 
By definition, $\LL_{X/Y}$ is $\WX$ equipped with the differential
$d_1$.  Since $d_2:\WX\to\S^2(\WX)[1]$ is a chain map of complexes we
get a natural map in the derived category
\[\LL_{X/Y}\to\S^2(\LL_{X/Y})[1],\]
which makes $\LL_{X/Y}[1]$ into a Lie coalgebra object in the derived
category of $\cO_X$-modules. 

\paragraph
Let us now return to the case of a lci embedding of schemes
$i:X\hookrightarrow Y$. In this case, $\LL_{X/Y}=N^\vee_{X/Y}[1]=
E[1]$.

\begin{Proposition}
\label{pp:HKRviaLinfty} 
The bracket
\[N^\vee_{X/Y}\to\wedge^2N^\vee_{X/Y}[2]\]
of Example~\ref{ex:bracket} is given by the HKR class $\alpha$.
\end{Proposition}

\begin{Proof}
As we saw above the complex $(\S(\WX),d)$ is filtered by complexes
$(\S^{\ge k}(\WX),d)$. As a graded sheaf the quotient $\S(\WX)^{\ge
  1}/\S(\WX)^{\ge 3}$ is equal to $\WX\bigoplus\S^2(\WX)$, however,
the differential on the quotient includes both the $d_1$ and the $d_2$
components. This provides an identification of the quotient
$\S(\WX)^{\ge 1}/\S(\WX)^{\ge 3}$ with the cone of the morphism
\[d_2:(\WX,d_1)\to(\S^2(\WX),d_1).\]

In our case, $(\S^k(\WX),d_1)$ has cohomology only in degree
$-k$. Therefore in the derived category $\D(X)$ we have
\[(\S^{\ge k}(\WX),d)\simeq\tau^{\le
  -k}((\S^k(\WX),d_1)\simeq\tau^{\le
  -k}(\cO_X\otimes_{i^{-1}\cO_Y}\cO_X)=\tau^{\le -k}(i^*i_*\cO_X).\]
Now the statement follows from Proposition~\ref{pp:alphadelta}.
\qed
\end{Proof}

\paragraph
\textbf{Remark. }  
There is a similar description of the HKR class $\alpha_V$ for any
vector bundle $V$ on $X$ (cf.\ the definition of the Atiyah classes in
\cite[Section~4.2]{BucFle}.  Namely we can resolve $V$ by a
$(\S(\W),d)$-module of the form
\[(\S(\W)\otimes_{\phi^{-1}\cO_Y}\RV,d).\]
Here $\RV$ is a graded sheaf of flat $\phi^{-1}\cO_Y$-modules.
The tensor product with $\cO_X$ gives a graded $\cO_X$-module 
\[(\S(\W)\otimes_{\phi^{-1}\cO_Y}\RV)\otimes_{\phi^{-1}\cO_Y}\cO_X=
\S(\WX)\otimes_{\phi^{-1}\cO_Y}\RV=\S(\WX)\otimes_{\cO_X}\RVX,\]
where $\RVX=\RV\otimes_{\phi^{-1}\cO_Y}\cO_X$. This tensor product
carries a differential $d$
compatible with the action of $(\S(\WX),d)$. 

This provides a coaction $d_k:\RVX\to\bwed^k(\WX[1])\otimes\RVX[1-k]$
of $\WX[1]$ on $\RVX$. In particular $d_0:\RVX\to\RVX[1]$ is a
differential, and it is easy to see that it turns $\RVX$ into a
resolution of $V$. Therefore $d_1$ yields a coaction of
$\LL_{X/Y}[1]$ on $V$ in $\D(X)$; one can check that the coaction is
given by the HKR class $\alpha_V$.  This is a (dual) analogue of the
usual statement that $\gog=T_X[-1]$ is a Lie algebra object in
$\D(X)$ and every object $F\in\D(X)$ is a representation of $\gog$.

\paragraph
\textbf{Remark.} 
Suppose that $i:X\hookrightarrow Y$ is an embedding of
quasi-projective schemes.  In this case, $\W$ can be chosen to be
locally free over $i^{-1}\cO_Y$, so that $(\WX,d_1)$ is a complex of
locally free $\cO_X$-modules. It is perhaps more natural to resolve
the $\cO_Y$-algebra $i_*\cO_X$ instead of the $i^{-1}\cO_Y$-algebra
$\cO_X$: this leads to a resolution of $X$ by a dg-scheme smooth over
$Y$ in the sense of~\cite{CioKap}.

\paragraph
\textbf{An example of a non-vanishing HKR class. } 
We shall now argue that the HKR class $\alpha$ is non-zero for the
embedding of $X=\p1\times\p1$ into $Y=\p5$ using the very
ample line bundle $\cO(1)\boxtimes\cO(2)$ on $X$, thus showing that
our theory is non-trivial.

\begin{Lemma}
\label{lm:ex_eta}
Consider  the cotangent bundle short exact sequence
\[ 0 \ra E = N^\chk \ra \Omega_Y^1|_X \ra \Omega_X^1 \ra 0, \]
and the associated long exact sequence
\[\dots\to H^1(X,i^*\Omega_Y^1)\to
H^1(X,\Omega_X^1)\overset{H^1(\eta)}\lra H^2(X,E)\to\dots.\]
Then we have 
\[ H^1(\eta) \neq 0. \]
\end{Lemma}
\vspace*{-6mm}

\begin{Proof} 
Clearly,
\[H^1(X,\Omega_X^1)=H^1(X,(\Omega_\p1\boxtimes\cO_\p1)\oplus(\cO_\p1\boxtimes\Omega_\p1))=\bbk^2,\]
so it suffices to show that $\dim H^1(X,i^*\Omega_Y^1)<2$. On $Y=\p5$,
consider the exact sequence
\[0\to \Omega^1_\p5\to \cO_\p5(-1)^6\to \cO_\p5\to 0.\] Applying
$i^*$, we obtain the sequence
\[0\to i^*\Omega_Y\to (\cO_\p1(-1)\boxtimes\cO_\p1(-2))^6\to\cO_X\to
0.\] Looking at the corresponding long exact sequence of cohomology
groups, we see that
\[\dim H^k(X,i^*\Omega_Y)=\begin{cases}0,&k\ne 1\\1,&k=1\end{cases},\]
as required.
\qed
\end{Proof}

\paragraph
Notice that for a line bundle $L$ on $X$ its Atiyah class $\At_L$
equals its first Chern class $c_1(L)\in H^1(X,\Omega^1_X)$. Therefore
\[\alpha_L=\eta\circ\At_L=H^1(\eta) (c_1(L))\in H^2(X,E), \]
where the latter is interpreted as the value of the map $H^1(\eta)$
evaluated on $c_1(L) \in H^1(X, \Omega_X^1)$.  It is clear that
$\alpha_L$ is additive in $L$.

\begin{Lemma} 
For the line bundle $L=\bwed^3 E$ we have $\alpha_L\ne 0$.
\label{lm:ex_alpha}
\end{Lemma}

\begin{Proof} 
Set
\[\ell=i^*\cO_\p5(1)=\cO_\p1(1)\boxtimes\cO_\p1(2).\] 
Since $\ell$ is a line bundle on $X$ that extends to $Y$,
$\alpha_\ell=0$. On the other hand,
\[L=i^*\omega_Y\otimes\omega_X^{-1}\simeq
i^*\cO_\p5(-6)\otimes(\cO_\p1(2)\boxtimes\cO_\p1(2))=
\cO_\p1(-4)\boxtimes\cO_\p1(-10).\] 
In particular, we see that $c_1(L)$ and $c_1(\ell)$ form a basis of
$H^1(X,\Omega^1_X)$. Now Lemma~\ref{lm:ex_eta} implies that
\[\alpha_L=H^1(\eta)(c_1(L))\ne 0,\]
since otherwise the map $H^1(\eta)$ would be identically zero.
\qed
\end{Proof}

\begin{Corollary} The class $\alpha=\alpha_{E}\in
H^2(X,E)$ is not zero.
\end{Corollary}

\begin{Proof}
Equivalently, we have to show that the vector bundle
$E$ does not extend to the first infinitesimal neighborhood
$X'\subset Y$. This is true because by Lemma~\ref{lm:ex_alpha} its
determinant $L$ does not extend to $X'$.
\qed
\end{Proof}

\paragraph
We conclude this section by speculating how Theorem~\ref{thm:mainthm}
should be understood from the point of view of the
$L_\infty$-coalgebra structure on the cotangent complex.

Recall from~(\ref{subsec:linfty}) that the construction of the
cotangent complex endowed the explicit complex $\WX[1]$ with an
$L_\infty$-coalgebra structure.  This is encoded either as a
sequence of maps 
\[ d_k:\WX[1] \ra \wedge^k(\WX[1])[2-k] \] 
for $k\geq 1$, or as a degree one differential $d$ on $\S(\WX)$.  The
complex $(\S(\WX), d)$ represents the derived tensor product $\cO_X
\otimes_{\phi^{-1}\cO_Y} \cO_X$ or, in the case $\phi$ is an affine
morphism of schemes, the object $\phi^*\phi_* \cO_X$ of $\D(X)$.

\paragraph
Now consider the case of a closed embedding $i$ of an lci subscheme
$X$ of a smooth scheme $Y$.  Then the cohomology of the complex
$(\WX[1], d_1)$ is concentrated in degree $-2$, where it equals $E$.
Assume that we can overcome the technical difficulties associated with
applying the homological perturbation lemma to the coalgebra $\WX[1]$.
(Since we work in the category of sheaves on $X$ which has positive
homological dimension the usual theory can not be applied directly.)
Then by transfer of structure we get an $L_\infty$-coalgebra structure
on $E[2]$ which is encoded by structure maps (after adjusting degrees)
\[ d_k':E \ra \wedge^k E[k]. \] 
These maps can further be assembled to yield a square-zero, degree one
``differential'' in $\D(X)$ 
\[ d':\S(E[1]) \ra \S(E[1])[1]. \] 
The main feature of homological perturbation theory is that 
$(\S(E[1]), d')$ shall be homotopic to $(\WX[1], d)$, so in particular
they will represent the same object of $\D(X)$ which we saw earlier
is $i^*i_* \cO_X$.  (Note however that while the map $d'$
has all the classical properties of a differential, it is not an actual
chain map, so the resulting object is not a complex and 
the notion of quasi-isomorphism does not make sense without further
developments of homological algebra.)
\medskip

\noindent
The above considerations combined with Theorem~\ref{thm:mainthm}
motivate us to state the following conjecture.

\begin{Conjecture}
\label{conj:mainconj}
The $L_\infty$-coalgebra $E[2]$ has the property that if the classical
co-bracket $d_2' = \alpha_E$ vanishes, then all the higher co-brackets
$d_k'$ vanish for $k\geq 2$ as well.
\end{Conjecture}

\paragraph
Note that in the spectral sequence of Corollary~\ref{cor:maincor}, if
the differentials up to the $(k-1)$-st page vanish, the
differentials at the $k$-th page are maps 
\[ H^p(X, \wedge^q E^\chk) \ra H^{p+k}(X, \wedge^{q-k+1} E^\chk). \]
It seems reasonable to guess that these maps are contraction with
$d_k'$, with $d_k'$ regarded as an element of $H^k(X, \wedge^kE
\otimes E)$.  

If we believe this then Conjecture~\ref{conj:mainconj} immediately
implies Corollary~\ref{cor:maincor}.  It also implies
Theorem~\ref{thm:mainthm} as well: $\alpha_E = 0$ implies $d_2' =0$,
hence by the Conjecture $d_k'=0$ for all $k\geq 2$, so $d' = 0$. Since
the complex $i^*i_* \cO_X = (\S(\WX), d)$ is quasi-isomorphic to
\[ (\S(E[1]), d) = (\S(E[1]), 0) = \S(E[1]), \]
Theorem~\ref{thm:mainthm} follows.